# Diophantine approximation by negative continued fraction

Hiroaki ITO

### Abstract

We show that the growth rate of denominator $Q_n$ of the $n$-th convergent of negative expansion of $x$ and the rate of approximation:

$$\frac{\log n}{n} \log \left| x - \frac{P_n}{Q_n} \right| \to -\frac{\pi^2}{3} \quad \text{in measure.}$$

for a.e. $x$. In the course of the proof, we reprove known inspiring results that arithmetic mean of digits of negative continued fraction converges to 3 in measure, although the limit inferior is 2, and the limit superior is infinite almost everywhere.

## 1 Introduction

Let $N$ be a fixed non-zero integer. An expression

$$[a_0; a_1, a_2, a_3, \cdots]_N := a_0 + \cfrac{N}{a_1 + \cfrac{N}{a_2 + \cfrac{N}{a_3 + \cdots}}}$$

is called an $N$-continued fraction and the $a_i$ are called partial quotients or digits. One can show that if we assume $a_i \geq N$ when $N \geq 1$ (resp. $a_i \geq |N| + 1$ when $N \leq -1$), then the expression $x = [a_0; a_1, a_2, \cdots]_N$ is unique for an irrational $x$. It is retrieved as a coding of the map $T_N : [0, 1) \to [0, 1)$ defined by

$$T_N(x) := \frac{N}{x} - \left\lfloor \frac{N}{x} \right\rfloor, \quad x \neq 0; \quad T_N(0) = 0.$$

More precisely, putting $a_0 = 0$ if $N \geq 1$; $a_0 = 1$ if $N \leq -1$ and

$$a_i = \left\lfloor \frac{N}{T_N^{i-1}(x)} \right\rfloor, N \geq 1; \quad a_i = -\left\lfloor \frac{N}{T_N^{i-1}(1-x)} \right\rfloor, N \leq -1, \quad \text{for } i = 1, 2, \cdots,$$

we have $x = [a_0; a_1, a_2, \cdots]_N$ for all $x \in [0, 1) \setminus \mathbb{Q}$. It is known that $([0, 1), \mathcal{B}, \mu_N, T_N)$ is ergodic with respect to the absolutely continuous invariant measure:

$$\mu_N(A) = \begin{cases} \dfrac{1}{\log \frac{N+1}{N}} \displaystyle\int_A \frac{dx}{N+x}, & \text{if } N \in \mathbb{Z} \setminus \{0, -1\}, \\ \displaystyle\int_A \frac{dx}{1-x}, & \text{if } N = -1. \end{cases}$$

In particular when $N = 1$, Birkhoff's ergodic theorem implies:

**Theorem 1.1.** *For a.e. $x \in [0, 1)$ with $x = [0; a_1, a_2, \cdots]_1$, one has*

$$\lim_{n \to \infty} \frac{1}{n} \#\{1 \leq i \leq n; a_i = k\} = \frac{1}{\log 2} \log \left( 1 + \frac{1}{k(k+2)} \right) \quad \text{for each } k \in \mathbb{N}, \tag{1}$$

$$\lim_{n\to\infty} \frac{n}{\frac{1}{a_1} + \frac{1}{a_2} + \cdots + \frac{1}{a_n}} = 1.74\cdots, \tag{2}$$

$$\lim_{n\to\infty} \sqrt[n]{a_1 a_2 \cdots a_n} = 2.68\cdots, \tag{3}$$

$$\lim_{n\to\infty} \frac{a_1 + a_2 + \cdots + a_n}{n} = \infty, \tag{4}$$

$$\lim_{n\to\infty} \frac{1}{n} \log q_n = \frac{\pi^2}{12 \log 2}, \quad \lim_{n\to\infty} \frac{1}{n} \log \left| x - \frac{p_n}{q_n} \right| = -\frac{\pi^2}{6 \log 2}. \tag{5}$$

where $p_n/q_n = [0; a_1, a_2, \cdots, a_n]_1$.

Analogous formula are obtained for $N \neq -1$. However, since $\mu_{-1}([0,1)) = \infty$, Birkhoff's ergodic theorem does not apply when $N = -1$.

We show that formulas corresponding to Theorem 1.1. K. Dajani and C. Kraaikamp showed that for a.e. $x \in [0,1)$ with $x = [1; a'_1, a'_2, \cdots]_{-1}$,

$$\lim_{n\to\infty} \frac{n}{\frac{1}{a'_1} + \frac{1}{a'_2} + \cdots + \frac{1}{a'_n}} = 2, \quad \lim_{n\to\infty} \sqrt[n]{a'_1 a'_2 \cdots a'_n} = 2$$

in [6] by using a recoding formula

$$[0; a_1, a_2, \cdots]_1 = [1; \underbrace{2, \cdots, 2}_{(a_1 - 1) \text{ times}}, a_2 + 2, \cdots, \underbrace{2, \cdots, 2}_{(a_{2k-1} - 1) \text{ times}}, a_{2k} + 2, \cdots]_{-1}. \tag{6}$$

By the formula (6), we can see that for a.e. $x = [1; a'_1, a'_2, \cdots]_{-1} \in [0,1)$, one has

$$\lim_{n\to\infty} \frac{1}{n} \#\{1 \leq i \leq n; a'_i = 2\} = 1.$$

J. Aaronson showed that arithmetic mean converges to 3 in measure [4]. Also, J. Aaronson and H. Nakada showed that the limit inferior is 2 and the limit superior is infinite for almost everywhere [5]. We reprove those theorems by using the formula (6) and an idea in [1]. Also, we obtain a formula corresponding to (5), that is,

$$\frac{1}{n} \log Q_n \to 0 \quad \text{a.e.,} \quad \frac{1}{n} \log \left| x - \frac{P_n}{Q_n} \right| \to 0 \quad \text{a.e.}$$

where $P_n/Q_n = [1; a'_1, a'_2, \cdots]_{-1}$. Further, we show that

$$\frac{\log n}{n} \log Q_n \to \frac{\pi^2}{6} \quad \text{in measure,} \quad \frac{\log n}{n} \log \left| x - \frac{P_n}{Q_n} \right| \to -\frac{\pi^2}{3} \quad \text{in measure,}$$

and the limit inferior and the limit superior of $(\log n/n) \log Q_n$, $(\log n/n) \log |x - P_n/Q_n|$ have different value respectively for a.e. $x$.

## 2  Arithmetic mean of negative continued fraction digits

In this section, we provide more simple proof of the following result [4], [5].

**Proposition 2.1** (J. Aaronson, H. Nakada). *Let $[1; a'_1, a'_2, \cdots]_{-1}$ be a negative continued fraction expansion of $x \in [0,1)$. Then, for all $\varepsilon > 0$,*

$$\lim_{n\to\infty} \lambda(\{x \ : \ \left| \frac{1}{n} \sum_{k=1}^n a'_k - 3 \right| > \varepsilon\}) = 0$$

*i.e., the arithmetic mean converges to 3 in measure. But this doesn't converge to 3 a.e.. Moreover,*

$$\liminf_{n\to\infty} \frac{1}{n} \sum_{k=1}^n a'_k = 2, \quad \limsup_{n\to\infty} \frac{1}{n} \sum_{k=1}^n a'_k = \infty \quad a.e.$$

To prove this, we prepare some Lemma. We denote the Lebesgue measure by $\lambda$ in this paper, and we set

$$E\begin{pmatrix} k_1 & \cdots & k_n \\ r & \cdots & r_n \end{pmatrix} := \{x \in [0,1) : a_{k_1}(x) = r_1, \cdots, a_{k_n}(x) = r_n\}.$$

**Lemma 2.1.** *For all $r_1, \cdots, r_i, r \in \mathbb{N}$*

$$\frac{\lambda\left(E\begin{pmatrix} k_1 & \cdots & k_n & k_{n+1} \\ r_1 & \cdots & r_n & r \end{pmatrix}\right)}{\lambda\left(E\begin{pmatrix} k_1 & \cdots & k_n \\ r_1 & \cdots & r_n \end{pmatrix}\right)} \ll \frac{1}{r^2} \tag{7}$$

*where $k_1, \cdots, k_n, k_{n+1}$ are all different from each other. If $k_1 < \cdots < k_n < k_{n+1}$, then there exists two constants $B$, $\beta > 0$ such that*

$$\left| \frac{\lambda\left(E\begin{pmatrix} k_1 & \cdots & k_n & k_{n+1} \\ r_1 & \cdots & r_n & r \end{pmatrix}\right)}{\lambda\left(E\begin{pmatrix} k_1 & \cdots & k_n \\ r_1 & \cdots & r_n \end{pmatrix}\right)} - \frac{\log\left(1 + \frac{1}{r(r+2)}\right)}{\log 2} \right| < Be^{-\beta\sqrt{k_{n+1} - k_n}}. \tag{8}$$

*Further, for all $k$,*

$$\sum_{r > \phi(n)} \lambda\left(E\begin{pmatrix} k \\ r \end{pmatrix}\right) \ll \frac{1}{\phi(n)} \tag{9}$$

*where "$\ll, \gg$" are Vinogradov symbols.*

Note that Lemma 2.1 can be generalized to following assertion holds in the same way as the proof of Khintchin's Theorem [1].

**Theorem 2.1.** *Let $[0; a_1, a_2, \cdots]_1$ be a regular continued fraction expansion of $x \in [0,1)$ and let $\Lambda : \mathbb{N} \to \mathbb{N}$ be a strictly increasing function. Then, for all $\varepsilon > 0$,*

$$\lim_{n \to \infty} \lambda(\{x \in [0,1) : \left| \frac{1}{n \log n} \sum_{i=1}^{n} a_{\Lambda(i)} - \frac{1}{\log 2} \right| > \varepsilon\}) = 0.$$

*i.e., $\sum_{i=1}^{n} a_{\Lambda(i)}(x)/(n \log n)$ converges to $1/\log 2$ in measure.*

Let $N = a_1 + a_3 + \cdots + a_{2k-1} + j$, $0 \le j < a_{2k-1}$, then by (6), we have

$$\frac{1}{N} \sum_{n=1}^{N} a'_n = 2 + \frac{\sum_{i=1}^{k} a_{2i}}{\sum_{i=1}^{k} a_{2i-1} + j}. \tag{10}$$

Therefore, by Theorem 2.1 and the following Proposition, it is sufficient to show that

$$\frac{j}{k \log k} \to 0 \quad \text{in measure.}$$

**Proposition 2.2.** *Let $(X, \mathcal{B}, \mu)$ be a measure space. Let $a_n, b_n : X \to \mathbb{R}_{>0}$ be measurable functions, and let $a, b$ be positive constants. Then, if $a_n \to a$, $b_n \to b$ in measure, then one has*

(a) $a_n + b_n \to a + b$ *in measure,*

(b) $a_n b_n \to ab$ *in measure,*

(c) $a_n / b_n \to a/b$ *in measure.*

*Proof.* We will prove that only (c). We set for all $\varepsilon > 0$,

$$A_n = \left\{ x : |a_n - a| > \frac{\varepsilon}{2b} \right\}, \quad B_n = \left\{ x : |b_n - b| > \frac{\varepsilon}{2a} \right\}.$$

Then for all $0 < \varepsilon' < b$, we have

$$\mu\left(\left\{x : \left|\frac{a_n}{b_n} - \frac{a}{b}\right| > \frac{\varepsilon}{b(b-\varepsilon')}\right\}\right) = \mu\left(\left\{x : |ba_n - ab_n| > \frac{\varepsilon b_n}{b - \varepsilon'}\right\}\right)$$
$$\leq \mu(\{x : |ba_n - ab_n| > \varepsilon\} \cup \{x : |ba_n - ab_n| > \frac{\varepsilon b_n}{b - \varepsilon'} \text{ and } b_n < b - \varepsilon'\})$$
$$\leq \mu(\{x : b|a_n - a| + a|b_n - b| > \varepsilon\}) + \mu(\{x : |b - b_n| > \varepsilon'\})$$
$$\leq \mu(A_n \cup B_n) + \mu(\{x : |b - b_n| > \varepsilon'\}) \to 0$$

$\square$

Now let

$$E_k = \{x \mid \frac{a_{2k+1}}{k} > \varepsilon\}, \quad e_{i,k} = \{x \mid a_{2i+1} < k^2\}, \quad F_k = \bigcap_{i=1}^{k} e_{i,k},$$

then by Lemma 2.1 (7), (9) we have

$$\lambda(E_k \cap F_k^c) \leq \lambda(F_k^c) \leq \sum_{i=1}^{k} \lambda(e_{i,k}^c) \ll k \cdot \frac{1}{k^2} = \frac{1}{k},$$

$$\lambda(E_k \cap F_k) < \frac{1}{\varepsilon} \cdot \frac{1}{k} \int_{F_k} a_{2k+1} dx < \frac{1}{\varepsilon} \cdot \frac{1}{k} \sum_{r=1}^{k^2} r \cdot \lambda\left(E\binom{2k+1}{r}\right) \ll \frac{\log k^2}{\varepsilon k}.$$

Thus, $j/(k \log k) < a_{2k+1}/k \to 0$ in measure. This prove the first assertion of the Proposition 2.1.

We can check that the following claim holds in the same way as the proof of Borel-Bernstein Theorem [8].

**Proposition 2.3** (Borel-Bernstein). *For a function $\varphi : \mathbb{N} \to (0, \infty)$ we set*

$$\mathcal{W}_\varphi = \{x = [1; a_1', a_2', \cdots]_{-1} \in (0,1) : a_n' > \varphi(n) \text{ for infinitely many } n \in \mathbb{N}\}$$

*If the series $\sum_{n=1}^{\infty} 1/\varphi(n)$ diverges, then*

$$\lambda(\mathcal{W}_\varphi^c) = 0.$$

By Proposition 2.3 and (10), we have

$$\limsup_{k \to \infty} \frac{\sum_{i=1}^{k} a_{2i}}{\sum_{i=1}^{k} a_{2i-1}} = \infty \quad a.e.$$

Further, since $T_1(x) = [0; a_2, a_3, \cdots]_1$, we have

$$\frac{1}{N} \sum_{i=1}^{N} a_i'(T_1(x)) = 2 + \frac{\sum_{i=1}^{k} a_{2i+1}}{\sum_{i=1}^{k} a_{2i} + j}$$

where $N = a_2 + a_4 + \cdots + a_{2k} + j$, $0 \leq j < a_{2k+2}$, and hence

$$\liminf_{N \to \infty} \frac{\sum_{i=1}^{k} a_{2i}}{\sum_{i=1}^{k} a_{2i-1} + j} = \liminf_{k \to \infty} \frac{\sum_{i=1}^{k} a_{2i}}{\sum_{i=1}^{k} a_{2i-1} + a_{2k+1} - 1} \leq \liminf_{k \to \infty} \frac{\sum_{i=1}^{k} a_{2i}}{\sum_{i=1}^{k} a_{2i+1}} = 0 \quad a.e.$$

which finishes the proof of Proposition 2.1.

For $x = [1; a_1', \cdots, a_n', \cdots]_{-1} = [0; a_1, \cdots, a_n, \cdots]_1$, a generalized mean map $M_{p,n}(x)$ defined by

$$M_{p,n}(x) := \left(\frac{\sum_{k=1}^{n} a_k'^p}{n}\right)^{\frac{1}{p}}.$$

Let $p$ be a number in $(0,1)$. Then, for $N = a_1 + a_3 + \cdots + a_{2k-1} + j$, $0 \leq j < a_{2k+1}$, we have

$$\frac{1}{N}\sum_{1=1}^{N} a_i'^p = 2^p - \frac{k+j}{N} + \frac{\frac{1}{k}\sum_{i=1}^{k}(a_{2i}+2)^p}{\frac{1}{k}N} \to 2^p \quad \text{a.e.}$$

Therefore,

$$\left(\frac{1}{n}\sum_{1=1}^{n} a_i'^p\right)^{\frac{1}{p}} \to 2 \quad \text{a.e.}$$

Note that $\lim_{p\to 0} M_{p,n} = \sqrt[n]{a_1' a_2' \cdots a_n'}$, $\lim_{p\to\infty} M_{p,n} = \max\{a_i\}_{i=1,\cdots n}$, $\lim_{p\to -\infty} M_{p,n} = \min\{a_i\}_{i=1,\cdots n}$ and, $M_{p,n}(x) \leq M_{q,n}(x)$ for $p < q$. Thus, the following table is obtained:

|  | $p < 1$ | $p = 1$ | $1 < p < \infty$ | $p \to \infty$ |
|---|---|---|---|---|
| $\lim_{n\to\infty} M_{p,n}(x)$ | 2 | 3 | ? | $\infty$ |
| type of convergence | a.e. | in measure | ? | a.e. |

Let $p \geq 2$ be a integer. Then, for $N = a_1 + a_3 + \cdots + a_{2k-1} + j$, $0 \leq j < a_{2k+1}$, we have

$$\frac{1}{N}\sum_{k=1}^{N} a_k'^p = 2^p + \frac{\sum_{l=1}^{p} 2^{p-l}\binom{p}{l}\sum_{i=1}^{k} a_{2i}^l}{\sum_{i=1}^{k} a_{2i-1} + j}.$$

But its limit inferior cannot be obtained in the same way as in the case $p = 1$. Also we don't know if it converges in measure. In general, we can calculate generalized mean for digits of $N$-continued fraction expansion $x = [a_0; a_1, \cdots]_N$ when $N \neq -1$.

## 3  Asymmetric behavior of denominator of convergent

Let $p_n/q_n = [0; a_1, a_2, \cdots, a_n]_1$, $P_n/Q_n = [1; a_1', a_2', \cdots, a_n']_{-1}$.

**Lemma 3.1.** *For all $x \in [0,1)$ with $x = [0; a_1, a_2, \cdots]_1 = [1; a_1', a_2', \cdots]_{-1}$, let $N = a_1 + a_3 + \cdots + a_{2k-1} + j$, $0 \leq j < a_{2k+1}$, then we have*

$$Q_N = q_{2k-1} + (j+1)q_{2k}, \tag{11}$$

$$\frac{1}{Q_N Q_{N-1}} < \left|x - \frac{P_{N-1}}{Q_{N-1}}\right| < \frac{1}{Q_{N-1}(Q_N - Q_{N-1})}. \tag{12}$$

*Proof.* We can prove (11) by induction, and the inequality (12) by

$$\frac{P_N - P_{N-1}}{Q_N - Q_{N-1}} < x < \frac{P_N}{Q_N} < \frac{P_{N-1}}{Q_{N-1}}, \quad P_N Q_{N-1} - P_{N-1} Q_N = 1, \quad \text{for all } N$$

$\square$

**Theorem 3.1.** *For a.e. $x \in [0,1)$ with $x = [1; a_1', a_2', \cdots]_{-1}$,*

$$\lim_{n\to\infty} \frac{1}{n}\log Q_n = 0, \quad \lim_{n\to\infty} \frac{1}{n}\log\left|x - \frac{P_n}{Q_n}\right| = 0.$$

*Moreover,*

$$\frac{\log n}{n}\log Q_n \to \frac{\pi^2}{6} \quad \text{in measure}, \quad \frac{\log n}{n}\log\left|x - \frac{P_n}{Q_n}\right| \to -\frac{\pi^2}{3} \quad \text{in measure}. \tag{13}$$

*Proof.* Let $x \in [0,1)$ with $x = [0; a_1, a_2, \cdots]_1 = [1; a'_1, a'_2, \cdots]_{-1}$. For sufficiently large $N > 0$, we can write $N = a_1 + a_3 + \cdots + a_{2k-1} + j$ for $k \geq 1$, $0 \leq j < a_{2k+1}$. Then, we have

$$0 \leq \frac{1}{N} \log Q_N = \frac{1}{N} \log(q_{2k-1} + (j+1)q_{2k}) < \frac{1}{N}\{\log(j+2) + \log q_{2k}\}$$
$$< \frac{\log(a_{2k+1}+2)}{N} + \frac{2k}{N} \cdot \frac{1}{2k} \log q_{2k}$$
$$< \frac{\frac{1}{k}\log(a_{2i+1}+2)}{\frac{1}{k}\sum_{i=1}^{k} a_{2i-1}} + \frac{2}{\frac{1}{k}\sum_{i=1}^{k} a_{2i-1}} \cdot \frac{1}{2k} \log q_{2k} \to 0 \quad \text{a.e.}$$

And,

$$0 \geq \frac{1}{N} \log \left| x - \frac{P_N}{Q_N} \right| > -\frac{1}{N} \log Q_N Q_{N+1} = -\left( \frac{1}{N} \log Q_N + \frac{N+1}{N} \cdot \frac{1}{N+1} \log Q_{N+1} \right) \to 0 \quad \text{a.e.}$$

Also, we have

$$\frac{\log k}{N} \log Q_N = \frac{\log k}{N} \log(q_{2k-1} + (j+1)q_{2k}) < \frac{\log k}{N}\{\log(j+2) + \log q_{2k}\}$$
$$< \frac{\frac{1}{k}\log(a_{2i+1}+2)}{\frac{1}{k \log k}\sum_{i=1}^{k} a_{2i-1}} + \frac{2}{\frac{1}{k \log k}\sum_{i=1}^{k} a_{2i-1}} \cdot \frac{1}{2k} \log q_{2k} \to \frac{\pi^2}{6} \quad \text{in measure,} \quad (14)$$

$$\frac{\log k}{N} \log Q_N = \frac{\log k}{N} \log(q_{2k-1} + (j+1)q_{2k}) > \frac{\log k}{N} \log q_{2k-1}$$
$$> \frac{2}{\frac{1}{k \log k}\sum_{i=1}^{k} a_{2i-1} + \frac{j}{k \log k}} \cdot \frac{2k-1}{2k} \cdot \frac{1}{2k-1} \log q_{2k-1} \to \frac{\pi^2}{6} \quad \text{in measure.}$$

Since if $a_n < b_n \to b$, $a_n > b'_n \to b$ in measure, then $a_n \to b$ in measure,

$$\frac{\log k}{N} \log Q_N \to \frac{\pi^2}{6} \quad \text{in measure.}$$

And, since

$$1 \leq \frac{\log N}{\log k} = \frac{\log \frac{N \log 2}{k \log k} + \log \frac{k \log k}{\log 2}}{\log k} \to 1 \quad \text{in measure,} \quad (15)$$

by Proposition 2.2,

$$\frac{\log n}{n} \log Q_n \to \frac{\pi^2}{6} \quad \text{in measure.}$$

Also, by Lemma 3.1,

$$\frac{\log N}{N} \log \left| x - \frac{P_N}{Q_N} \right| > -\frac{\log N}{N} \log Q_N - \frac{\log N}{N} \log Q_{N+1} \to -\frac{\pi^2}{3} \quad \text{in measure,}$$

$$\frac{\log N}{N} \log \left| x - \frac{P_N}{Q_N} \right| < -\frac{\log N}{N} \log Q_N - \frac{2\frac{\log N}{\log k}}{\frac{1}{k \log k} N} \frac{1}{2k} \log q_{2k} \to -\frac{\pi^2}{3} \quad \text{in measure.}$$

Thus,

$$\frac{\log n}{n} \log \left| x - \frac{P_n}{Q_n} \right| \to -\frac{\pi^2}{3} \quad \text{in measure.}$$

□

To obtain limit inferior and limit superior of (13), we use the following Lemma (see [3]).

**Lemma 3.2** (Diamond, Vaaler). *For almost all* $x = [0, a_1, a_2, \cdots ]_1 \in (0,1)$,

$$\lim_{n \to \infty} \frac{1}{n \log n} \left( \sum_{i=1}^{n} a_{2i-1} - \max_{1 \leq i \leq n} a_{2i-1} \right) = \frac{1}{\log 2}.$$

**Theorem 3.2.**

$$\liminf_{n \to \infty} \frac{\log n}{n} \log Q_n = 0, \quad \limsup_{n \to \infty} \frac{\log n}{n} \log \left| x - \frac{P_n}{Q_n} \right| = 0 \quad a.e.$$

*Proof.* By Borel-Bernstein theorem and (14),

$$\liminf_{n \to \infty} \frac{\log k}{N} \log Q_N = 0 \quad \text{a.e.}$$

Also, since there are at most finitely many $k$ such that $a_k > k^2$ by Borel-Bernstein theorem, for sufficiently large $k$,

$$\log \left( \frac{\log 2 \cdot N}{k \log k} \right) < \log \left( \frac{\log 2 \cdot \left( \sum_{i=1}^{k} a_{2i-1} - \max_{1 \leq i \leq k} \{a_{2i-1}\} \right)}{k \log k} + \frac{\log 2 \cdot a_{2k+1}}{k \log k} + \frac{\log 2 \cdot \max_{1 \leq i \leq k} \{a_{2i+1}\}}{k \log k} \right)$$

$$< \log \left( \frac{\log 2 \cdot \left( \sum_{i=1}^{k} a_{2i-1} - \max_{1 \leq i \leq k} \{a_{2i-1}\} \right)}{k \log k} + \frac{\log 2 \cdot k}{\log k} + \frac{\log 2 \cdot k}{\log k} \right),$$

Therefore, by Lemma 3.2 and (15),

$$\lim_{k \to \infty} \frac{\log N}{\log k} = 1 \quad \text{a.e.} \tag{16}$$

Thus, the claim is proved. $\square$

**Proposition 3.1** (Riesz). *Let $(X, \mathcal{B}, \mu)$ be a measure space. If $f_n, f : X \to \mathbb{R}$ are measurable function such that $f_n \to f$ in measure, then there exist subsequence $(n_k)_{k \in \mathbb{N}}$ such that $f_{n_k} \to f$ a.e.*

**Lemma 3.3.**

$$\liminf_{k \to \infty} \frac{1}{k \log k} \sum_{i=1}^{k} a_{2i-1} = \frac{1}{\log 2} \quad a.e.$$

*Proof.* By Lemma 3.2,

$$\frac{1}{k \log k} \sum_{i=1}^{k} a_{2i-1} > \frac{1}{k \log k} \left( \sum_{i=1}^{k} a_{2i-1} - \max_{1 \leq i \leq k} a_{2i-1} \right) \to \frac{1}{\log 2} \quad \text{a.e.}$$

On the other hand, by Proposition 3.1,

$$\liminf_{k \to \infty} \frac{1}{k \log k} \sum_{i=1}^{k} a_{2i-1} \leq \frac{1}{\log 2} \quad \text{a.e.,}$$

which finishes the proof of this Lemma. $\square$

**Theorem 3.3.**

$$\limsup_{n \to \infty} \frac{\log n}{n} \log Q_n = \frac{\pi^2}{6} \quad a.e., \quad \liminf_{n \to \infty} \frac{\log n}{n} \log \left| x - \frac{P_n}{Q_n} \right| = -\frac{\pi^2}{3} \quad a.e.$$

*Proof.* By Borel-Bernstein theorem, for sufficiently large $k$,

$$\frac{1}{k} \log (a_{2k+1} + 2) < \frac{1}{k} \log ((2k+1)^2 + 2),$$

and then by (14), Lemma 3.3 and Proposition 3.1,

$$\limsup_{N\to\infty} \frac{\log k}{N} \log Q_N = \frac{\pi^2}{6} \quad \text{a.e.}$$

Thus, by (16),

$$\limsup_{n\to\infty} \frac{\log n}{n} \log Q_n = \frac{\pi^2}{6} \quad \text{a.e.}$$

Also, by Lemma 3.1, we have

$$\liminf_{n\to\infty} \frac{\log n}{n} \log \left| x - \frac{P_n}{Q_n} \right| = -\frac{\pi^2}{3} \quad \text{a.e.}$$

□

Note that

$$\prod_{i=0}^{n-1} T_{-1}^i(x) = |Q_{n-1}(1-x) - P_{n-1}|$$

(where $P_n/Q_n$ is $n$th-convergent of negative expansion of $1-x$), then by (13), we have

$$\frac{\log n}{n} \sum_{i=0}^{n-1} \log T_{-1}^i(x) \to \int_{[0,1)} \log x \, d\mu_{-1}(x) \quad \text{in measure.}$$